\documentclass{svmult}

\usepackage{amsmath,amssymb,latexsym,amsfonts}

\newtheorem{Thm}{Theorem}[section]

\newtheorem{Prop}[Thm]{Proposition}
\newtheorem{Cor}[Thm]{Corollary}

\newtheorem{Rem}[Thm]{Remark}
\newtheorem{Ex}[Thm]{Example}

\newcommand{\eps}{\ensuremath{\varepsilon}}
\newcommand{\e}{{\rm e}}
\renewcommand{\d}{{\rm d}}

\newcommand{\law}{\stackrel{{\rm law}}{=}}
\newcommand{\claw}{\stackrel{{\rm law}}{\longrightarrow}}
\newcommand{\cdist}{\stackrel{{\rm d}}{\longrightarrow}}

\newcommand{\bR}{\ensuremath{\mathbb{R}}}
\newcommand{\bZ}{\ensuremath{\mathbb{Z}}}

\newcommand{\cB}{\ensuremath{\mathcal{B}}}
\newcommand{\cM}{\ensuremath{\mathcal{M}}}

\newcommand{\absol}[1]{\left| #1 \right|} 
\newcommand{\rbra}[1]{\left( #1 \right)} 
\newcommand{\cbra}[1]{\left\{ #1 \right\}} 
\newcommand{\sbra}[1]{\left[ #1 \right]} 

\begin{document}

\title*{On a zero-one law for the norm process of transient random walk}
\author{Ayako Matsumoto\inst{1}\and
Kouji Yano\inst{1}}
\institute{T \& D Financial Life Insurance Company, JAPAN. 
\and Graduate School of Science, Kobe University, Kobe, JAPAN.}

\maketitle

\renewcommand{\theequation}{\arabic{section}.\arabic{equation}}
\newcounter{bean}


\begin{abstract}
A zero-one law of Engelbert--Schmidt type 
is proven for the norm process of a transient random walk. 
An invariance principle for random walk local times 
and a limit version of Jeulin's lemma play key roles. 
\end{abstract}


\section{Introduction}

Let $ S=(S_n:n \in \bZ_{\ge 0}) $ be a random walk in $ \bZ^d $ 
starting from the origin. 
Let $ \| \cdot \| $ be a norm on $ \bR^d $ 
taking integer values on the integer lattice $ \bZ^d $. 
The norm $ \| \cdot \| $ 
cannot be the Euclidean norm 
denoted by $ |x| = \sqrt{|x^1|^2 + \cdots + |x^d|^2} $. 
By the {\em norm process of the random walk} $ S $, 
we mean the process $ \| S \| = (\| S_n \|: n \in \bZ_{\ge 0}) $. 
The purpose of the present paper is to study summability of $ f(\| S_n \|) $ 
for a non-negative function $ f $ on $ \bZ $. 

Set $ X_n=S_n-S_{n-1} $ for $ n \in \bZ_{\ge 1} $. 
Then $ X_n $'s 
are independent identically-distributed random vectors taking values in $ \bZ^d $. 
We suppose that $ E[X_1^i] = 0 $ and $ E[(X_1^i)^2] < \infty $, $ i=1,2,\ldots,d $. 
Let $ Q $ denote the covariance matrix of $ X_1 $, 
i.e., $ Q = (E[X_1^i X_1^j])_{i,j} $. 
We introduce the following assumption: 
\\ \quad {\rm (A0)} 
$ Q = \sigma^2 I $ for some constant $ \sigma>0 $, 
where $ I $ stands for the identity matrix. 

We write 
\begin{align}
B(0;r) =& \{ x \in \bR^d : \| x \| \le r \} , 
\label{} \\
\partial B(0;r) =& \{ x \in \bR^d : \| x \| = r \} . 
\label{}
\end{align}
For $ k \in \bZ_{\ge 0} $, we set 
\begin{align}
N(k) = \sharp (\partial B(0;k) \cap \bZ^d) 
= \sharp \cbra{ x \in \bZ^d : \| x \| = k } . 
\label{}
\end{align}
We call $ B $ a {\em $ d $-polytope} 
if $ B $ is a bounded convex region in a $ d $-dimensional space 
enclosed by a finite number of $ (d-1) $-dimensional hyperplanes. 
The part of the polytope $ B $ which lies in one of the hyperplanes 
is called a {\em cell}. 
(See, e.g., \cite{Cox} for this terminology.) 
We introduce the following assumptions: 
\\ \quad {\rm (A1)} 
$ \| x \| \in \bZ_{\ge 0} $ for any $ x \in \bZ^d $. 
\\ \quad {\rm (A2)} 
For each $ k \in \bZ_{\ge 1} $, 
the set $ B(0;k) $ is a $ d $-polytope 
whose vertices are contained in $ \bZ^d $. 
Consequently, its boundary $ \partial B(0;k) $ 
is the union of all cells of the $ d $-polytope $ B(0;k) $. 
\\ \quad {\rm (A3)} 
For any $ k \in \bZ_{\ge 1} $, 
there exists a finite partition of $ \partial B(0;1) $, 
which is denoted by $ \{ U^{(k)}_j:j=1,\ldots,M(k) \} $, 
such that the following statements hold: 
\\ \quad \quad {\rm (i)} 
$ M(k) \le N(k) $ and $ M(k)/N(k) \to 1 $ as $ k \to \infty $; 
\\ \quad \quad {\rm (ii)} 
each $ U^{(k)}_j $ contains at least one point of $ \partial B(0;1) \cap (k^{-1} \bZ^d) $; 
\\ \quad \quad {\rm (iii)} 
the $ U^{(k)}_j $'s for $ j=1,\ldots,M(k) $ have a common area; 
\\ \quad \quad {\rm (iv)} 
$ \max_{j} \max \{ \| x-y \| : x,y \in U^{(k)}_j \} \to 0 $ as $ k \to \infty $. 

Note that these assumptions {\rm (A0)}-{\rm (A3)} imply that 
$ N(k) \to \infty $ as $ k \to \infty $. 
Our main theorem is the following: 

\begin{Thm} \label{thm: main2} 
Suppose that $ d \ge 3 $ and that {\rm (A0)}-{\rm (A3)} hold. 
Then, for any non-negative function $ f $ on $ \bZ_{\ge 0} $, 
the following conditions are equivalent: 
\\ \quad {\rm (I)} 
$ \displaystyle P \rbra{ \sum_{n=1}^{\infty } f(\| S_n \|) < \infty } > 0 $; 
\\ \quad {\rm (II)} 
$ \displaystyle P \rbra{ \sum_{n=1}^{\infty } f(\| S_n \|) < \infty } = 1 $; 
\\ \quad {\rm (III)} 
$ \displaystyle E \sbra{ \sum_{n=1}^{\infty } f(\| S_n \|) } < \infty $; 
\\ \quad {\rm (IV)} 
$ \displaystyle \sum_{k=1}^{\infty } k^{2-d} N(k) f(k) < \infty $. 
\\
Suppose, moreover, that 
\\ \quad {\rm (A4)} 
There exists $ k_0 \in \bZ_{\ge 1} $ such that $ N(k) $ is non-decreasing in $ k \ge k_0 $. 
\\
Then the above conditions are equivalent to 
\\ \quad {\rm (V)} 
$ \displaystyle \sum_{k=1}^{\infty } k f(k) < \infty $. 
\end{Thm}

We will prove, in Section \ref{sec: proof1}, 
that {\rm (III)} and {\rm (IV)} are equivalent, 
by virtue of the asymptotic behavior of the Green function 
due to Spitzer \cite{Spitzer} (see Theorem \ref{thm: Spitzer}). 
We will prove, in Section \ref{sec: proof2}, 
that {\rm (I)} implies {\rm (IV)}, 
where a key role is played by 
a limit version of {\em Jeulin's lemma} (see Proposition \ref{thm: limit Jeulin1}). 
Note that {\rm (III)} trivially implies {\rm (II)} 
and that {\rm (II)} trivially implies {\rm (I)}. 

\begin{center}
\begin{tabular}{rcccl}
        &                & \S \ref{sec: proof1}    &              & \\
        & {\rm (III)}     & $ \Longleftrightarrow $ & {\rm (IV)}   & \\
trivial & $ \Downarrow $ &                         & $ \Uparrow $ & \S \ref{sec: proof2} \\
        & {\rm (II)}     & $ \Longrightarrow $     & {\rm (I)}   & \\
        &                & trivial                 &              & 
\end{tabular}
\end{center}

The equivalence between {\rm (I)} and {\rm (II)} 
may be considered to be a {\em zero-one law of Engelbert--Schmidt type}; 
see Section \ref{sec: zero-one}. 
However, we remark that this equivalence follows 
also from the {\em Hewitt--Savage zero-one law} 
(see, e.g., \cite[Thm.7.36.5]{Billingsley}). 
In fact, the event $ \{ \sum f(\| S_n \|) < \infty \} $ is {\em exchangeable}, 
i.e., invariant under permutation of any finite number of the sequence $ (X_n) $.

If $ d=1 $ or $ 2 $, the random walk $ S $ is recurrent, 
and hence it is obvious that 
the conditions {\rm (I)}-{\rm (III)} are equivalent to 
stating that $ f(k) \equiv 0 $. 
This is why we confine ourselves to the case $ d \ge 3 $, 
where the random walk $ S $ is transient 
so that $ \| S_n \| $ diverges as $ n \to \infty $. 
In the case $ d \ge 3 $, the summability of $ f(\| S_n \|) $ depends upon 
how rapidly the function $ f(k) $ vanishes as $ k \to \infty $. 
Theorem \ref{thm: main2} gives a criterion 
for the summability of $ f(\| S_n \|) $ 
in terms of summability of $ kf(k) $.

Consider the {\em max norm} 
\begin{align}
\| x \|^{(d)}_{\infty } = \max_{i=1,\ldots,d} |x^i| 
, \quad x = (x^1,\ldots,x^d) \in \bR^d 
\label{}
\end{align}
and the {\em $ \ell^1 $-norm} 
\begin{align}
\| x \|^{(d)}_1 = \sum_{i=1}^d |x^i| 
, \quad x = (x^1,\ldots,x^d) \in \bR^d . 
\label{}
\end{align}
We will show in Section \ref{sec: norm} that 
these norms satisfy {\rm (A1)}-{\rm (A4)}. 
Thus we obtain the following corollary: 

\begin{Cor} \label{thm: main1}
Let $ S $ be a simple random walk of dimension $ d \ge 3 $ 
and take $ \| \cdot \| $ as the max norm or the $ \ell^1 $-norm. 
Then, for any non-negative function $ f $ on $ \bZ_{\ge 0} $, 
the conditions {\rm (I)}-{\rm (V)} are equivalent. 
\end{Cor}

The organization of this paper is as follows. 
In Section \ref{sec: zero-one}, 
we give a brief summary of known results 
of zero-one laws of Engelbert--Schmidt type. 
In Section \ref{sec: Jeulin}, 
we recall Jeulin's lemma. 
We also state and prove its limit version in discrete time. 
In Section \ref{sec: norm}, 
we present some examples of norms which satisfy {\rm (A1)}-{\rm (A4)}. 
Sections \ref{sec: proof1} and \ref{sec: proof2} 
are devoted to the proof of Theorem \ref{thm: main2}. 
In Section \ref{sec: Shiga}, 
we present some results about Jeulin's lemma 
obtained by Shiga \cite{Shiga}.

\section{Zero-one laws of Engelbert--Schmidt type} \label{sec: zero-one}
\setcounter{equation}{0}

Let us give a brief summary of known results of zero-one laws 
concerning finiteness of certain integrals, 
which we call {\em zero-one laws of Engelbert--Schmidt type}. 

\

\noindent
$ 1^{\circ}) $. 
Let $ (B_t:t \ge 0) $ be a one-dimensional Brownian motion starting from the origin. 
The following theorem, 
which originates from Shepp--Klauder--Ezawa \cite{SKE} 
with motivation in quantum theory, 
is due to Engelbert--Schmidt \cite[Thm.1]{ES1} 
with motivation in construction of a weak solution 
of a certain stochastic differential equation 
by means of time-change method. 

\begin{Thm}[\cite{SKE},\cite{ES1}] \label{thm: ES}
Let $ f $ be a non-negative Borel function on $ \bR $. 
Then the following conditions are equivalent: 
\\ \quad {\rm (B1)} 
$ P \rbra{ \int_0^t f(B_s) \d s < \infty \ \text{for every $ t \ge 0 $} } > 0 $; $ \Big. $
\\ \quad {\rm (B2)} 
$ P \rbra{ \int_0^t f(B_s) \d s < \infty \ \text{for every $ t \ge 0 $} } = 1 $; $ \Big. $
\\ \quad {\rm (B3)} 
$ f(x) $ is integrable on all compact subsets of $ \bR $. $ \Big. $
\end{Thm}

The proof of Theorem \ref{thm: ES} 
was based on the formula 
\begin{align}
\int_0^t f(B_s) \d s = \int_{\bR} f(x) L^B_t(x) \d x 
\label{}
\end{align}
where $ L^B_t(x) $ stands for the local time at level $ x $ by time $ t $ 
(see \cite{IM}). 

Engelbert--Schmidt \cite[Thm.1]{ES2} 
proved that a similar result holds for a Bessel process of dimension $ d \ge 2 $ 
starting from a positive number. 

\

\noindent
$ 2^{\circ}) $. 
Let $ (R_t:t \ge 0) $ be a Bessel process of dimension $ d > 0 $ 
starting from the origin, i.e., $ R_t = \sqrt{Z_t} $ 
where $ Z_t $ is the unique non-negative strong solution of 
\begin{align}
Z_t = t d + 2 \int_0^t \sqrt{|Z_s|} \d B_s . 
\label{}
\end{align}
The following theorem is due to 
Pitman--Yor \cite[Prop.1]{PY2} and Xue \cite[Prop.2]{Xue}. 

\begin{Thm}[\cite{PY2},\cite{Xue}] \label{thm: Xue0}
Suppose that $ d \ge 2 $. 
Let $ f $ be a non-negative Borel function on $ [0,\infty ) $. 
Then the following conditions are equivalent: 
\\ \quad {\rm (R1)} 
$ P \rbra{ \int_0^t f(R_s) \d s < \infty \ \text{for every $ t \ge 0 $} } > 0 $; $ \Big. $
\\ \quad {\rm (R2)} 
$ P \rbra{ \int_0^t f(R_s) \d s < \infty \ \text{for every $ t \ge 0 $} } = 1 $; $ \Big. $
\\ \quad {\rm (R3)} 
$ f(r) $ is integrable on all compact subsets of $ (0,\infty ) $ and $ \Big. $
\\ \quad \quad {\rm (R3a)} 
$ \int_0^c f(r) r (\log \frac{1}{r})_+ \d r < \infty $ if $ d=2 $; $ \Big. $
\\ \quad \quad {\rm (R3b)} 
$ \int_0^c f(r) r \d r < \infty $ if $ d>2 $ $ \Big. $ 
\\
where $ c $ is an arbitrary positive number. 
\end{Thm}

The proof of Theorem \ref{thm: Xue0} was done 
by applying Jeulin's lemma (see Theorem \ref{thm: Jeulin} below) 
to the total local time, 
where the assumption of Jeulin's lemma 
was assured by the {\em Ray--Knight theorem} (see Le Gall \cite[pp.299]{LeGall}). 

\

\noindent
$ 3^{\circ}) $. 
Xue \cite[Cor.4]{Xue} generalized Engelbert--Schmidt \cite[Cor. on pp.227]{ES2} 
and proved the following theorem. 

\begin{Thm}[\cite{Xue}] \label{thm: Xue}
Suppose that $ d>2 $ 
Let $ f $ be a non-negative Borel function on $ [0,\infty ) $. 
Then the following conditions are equivalent: 
\\ \quad {\rm (RI)} 
$ P \rbra{ \int_0^{\infty } f(R_t) \d t < \infty \Big. } > 0 $; 
\\ \quad {\rm (RII)} 
$ P \rbra{ \int_0^{\infty } f(R_t) \d t < \infty \Big. } = 1 $; 
\\ \quad {\rm (RIII)} 
$ E \sbra{ \int_0^{\infty } f(R_t) \d t \Big. } < \infty $; 
\\ \quad {\rm (RIV)} 
$ \int_0^{\infty } r f(r) \d r < \infty $. $ \Big. $
\end{Thm}

The proof of Theorem \ref{thm: Xue} was 
based on Jeulin's lemma and the Ray--Knight theorem. 
Our results (Theorem \ref{thm: main2} and Corollary \ref{thm: main1}) 
may be considered to be random walk versions of Theorem \ref{thm: Xue}. 
Note that, in Theorem \ref{thm: Xue}, the condition {\rm (RIII)}, 
which is obviously stronger than {\rm (RII)}, 
is in fact equivalent to {\rm (RII)}. 
We remark that, in Theorem \ref{thm: Xue}, 
we consider the perpetual integral $ \int_0^{\infty } f(R_t) \d t $ 
instead of the integrals on compact intervals.

\

\noindent
$ 4^{\circ}) $. 
H\"ohnle--Sturm \cite{HS1},\cite{HS2} 
obtained a zero-one law about the event 
\begin{align}
\mbox{$ \cbra{ \int_0^t f(X_s) \d s < \infty \ \text{for every $ t \ge 0 $} } $} 
\label{}
\end{align}
where $ (X_t:t \ge 0) $ is a symmetric Markov process 
which takes values in a Lusin space 
and which has a strictly positive density. 
Their proof was based on excessive functions. 
As an application, they obtained the following theorem (\cite[pp.411]{HS2}). 

\begin{Thm}[\cite{HS2}]
Suppose that $ 0<d<2 $. 
Let $ f $ be a non-negative Borel function on $ [0,\infty ) $. 
Then the following conditions are equivalent: 
\\ \quad {\rm (Ri)} 
$ P \rbra{ \int_0^t f(R_s) \d s < \infty \ \text{for every $ t \ge 0 $} } > 0 $; $ \Big. $
\\ \quad {\rm (Rii)} 
$ P \rbra{ \int_0^t f(R_s) \d s < \infty \ \text{for every $ t \ge 0 $} } = 1 $; $ \Big. $
\\ \quad {\rm (Riii)} 
$ f(x) $ is integrable on all compact subsets of $ [0,\infty ) $ 
and 
$ \int_0^1 f(x) x^{d-1} \d x < \infty $. $ \Big. $
\end{Thm}

See also Cherny \cite[Cor.2.1]{Cherny} for another approach.

\

\noindent
$ 5^{\circ}) $. 
Engelbert--Senf \cite{ES3} 
studied integrability of $ \int_0^{\infty } f(Y_s) \d s $ 
where $ (Y_t:t \ge 0) $ is a Brownian motion with constant drift. 
See Salminen--Yor \cite{SY1} for a generalization of this direction. 
See also Khoshnevisan--Salminen--Yor \cite{KSY} 
for a generalization of the case 
where $ (Y_t:t \ge 0) $ is a certain one-dimensional diffusion process.

\section{Jeulin's lemma and its limit version in discrete time} \label{sec: Jeulin}
\setcounter{equation}{0}

\subsection{Jeulin's lemma}

Jeulin \cite[Lem.3.22]{Jeulin1} gave quite a general theorem 
about integrability of a function of a stochastic process. 
He gave detailed discussions in \cite{Jeulin2} about his lemma. 
Among the applications presented in \cite{Jeulin2}, 
let us focus on the following theorem:

\begin{Thm}[\cite{Jeulin1},\cite{Jeulin2}] \label{thm: Jeulin}
Let $ (X(t):0 < t \le 1) $ be a non-negative measurable process 
and $ \varphi $ a positive function on $ (0,1] $. 
Suppose that 
there exists a random variable $ X $ with 
\begin{align}
E[X]<\infty \quad \text{and} \quad P(X>0) = 1 
\label{}
\end{align}
such that 
\begin{align}
\frac{X(t)}{\varphi(t)} \law X 
\quad \text{holds for each fixed $ 0 < t \le 1 $}. 
\label{eq: equal in law}
\end{align}
Then, for any non-negative Borel measure $ \mu $ on $ (0,1] $, 
the following conditions are equivalent: 
\\ \quad {\rm (JI)} 
$ P \rbra{ \int_0^1 X(t) \mu(\d t) < \infty \Big. } > 0 $; 
\\ \quad {\rm (JII)} 
$ P \rbra{ \int_0^1 X(t) \mu(\d t) < \infty \Big. } = 1 $; 
\\ \quad {\rm (JIII)} 
$ E \sbra{ \int_0^1 X(t) \mu(\d t) \Big. } < \infty $; 
\\ \quad {\rm (JIV)} 
$ \int_0^1 \varphi(t) \mu(\d t) < \infty $. $ \Big. $
\end{Thm}

A good elementary proof of Theorem \ref{thm: Jeulin} 
can be found in Xue \cite[Lem.2]{Xue}. 

For several applications of Jeulin's lemma (Theorem \ref{thm: Jeulin}), 
see Yor \cite{Yor}, 
Pitman--Yor \cite{PY1}, \cite{PY2}, 
Xue \cite{Xue}, 
Peccati--Yor \cite{PeY1}, 
Funaki--Hariya--Yor \cite{FHY1}, \cite{FHY2}, 
and 
Fitzsimmons--Yano \cite{FY}. 

We cannot remove the assumption $ E[X]<\infty $ from Theorem \ref{thm: Jeulin}; 
see Proposition \ref{thm: Shiga5}.

\subsection{A limit version of Jeulin's lemma in discrete time}

For our purpose, 
we would like to replace 
the assumption \eqref{eq: equal in law} which requires identity in law 
by a weaker assumption which requires convergence in law. 
The following proposition plays a key role in our purpose 
(see also Corollary \ref{thm: Shiga2}). 

\begin{Prop} \label{thm: limit Jeulin1}
Let $ (V(k):k \in \bZ_{\ge 1}) $ be a non-negative measurable process 
and $ \Phi $ a positive function on $ \bZ_{\ge 1} $. 
Suppose that 
there exists a random variable $ X $ with 
\begin{align}
P(X>0) = 1 
\label{}
\end{align}
such that 
\begin{align}
\frac{V(k)}{\Phi(k)} \claw X 
\quad \text{as $ k \to \infty $}. 
\label{eq: conv in law}
\end{align}
Then, for any non-negative function $ f $ on $ \bZ_{\ge 1} $, 
it holds that 
\begin{align}
P \rbra{ \sum_{k=1}^{\infty } f(k) V(k) < \infty } > 0 
\quad \text{implies} \quad 
\sum_{k=1}^{\infty } f(k) \Phi(k) < \infty. 
\label{}
\end{align}
\end{Prop}

The following proof of Proposition \ref{thm: limit Jeulin1} 
is a slight modification of that of \cite[Lem.2]{Xue}. 

\begin{proof}
Suppose that $ P(\sum f(k) V(k) < \infty ) > 0 $. 
Then there exists a number $ C $ such that the event 
\begin{align}
B = \cbra{ \sum_{k=1}^{\infty } f(k) V(k) \le C } 
\label{}
\end{align}
has positive probability. 
Since $ P(X \le 0) = 0 $, 
there exists a positive number $ u_0 $ such that 
$ P(X \le u_0) < P(B)/4 $. 
By assumption \eqref{eq: conv in law}, we see that 
there exists $ u_1 $ with $ 0<u_1<u_0 $ such that 
\begin{align}
P(V(k)/\Phi(k) \le u_1) 
\stackrel{k \to \infty }{\longrightarrow} 
P(X \le u_1) < \frac{1}{4} P(B) . 
\label{}
\end{align}
Then, for some large number $ k_0 $, we have 
\begin{align}
P(V(k)/\Phi(k) \le u_1) \le \frac{1}{2} P(B) 
, \quad k \ge k_0 . 
\label{}
\end{align}
Now we obtain 
\begin{align}
C 
\ge& E \sbra{ 1_B \sum_{k=1}^{\infty } f(k) V(k) } 
\label{} \\
=& \sum_{k=1}^{\infty } f(k) \Phi(k) E \sbra{ 1_B \cdot \frac{V(k)}{\Phi(k)} } 
\label{} \\
=& \sum_{k=1}^{\infty } f(k) \Phi(k) \int_0^{\infty } P( B \cap \{ V(k)/\Phi(k) > u \} \d u 
\label{} \\
\ge& \sum_{k=k_0}^{\infty } f(k) \Phi(k) \int_0^{u_1} 
\sbra{ P(B) - P(V(k)/\Phi(k) \le u) }_+ \d u 
\label{} \\
\ge& \frac{1}{2} P(B) u_1 \sum_{k=k_0}^{\infty } f(k) \Phi(k) . 
\label{}
\end{align}
Since $ P(B) u_1 > 0 $, we conclude that $ \sum f(k) \Phi(k) < \infty $. 
\end{proof}

\section{Examples of norms} \label{sec: norm}
\setcounter{equation}{0}

Let us introduce several notations. 
For an index set $ A $ 
(we shall take $ A = \bZ_{\ge 0} $ or $ \bZ^d \setminus \{ 0 \} $ later), 
we denote $ \cM(A) $ 
the set of all non-negative functions on $ A $. 
For three functions $ f,g,h \in \cM(A) $ , 
we say that 
\begin{align}
f(a) \sim g(a) 
\quad \text{as $ h(a) \to \infty $} 
\label{}
\end{align}
if $ f(a)/g(a) \to 1 $ as $ h(a) \to \infty $. 
For two functions $ f,g \in \cM(A) $, 
we say that 
\begin{align}
f(a) \asymp g(a) 
\quad \text{for $ a \in A $} 
\label{}
\end{align}
if there exist positive constants $ c_1,c_2 $ such that 
\begin{align}
c_1 f(a) \le g(a) \le c_2 f(a) 
\quad \text{for $ a \in A $} . 
\label{}
\end{align}
For two functionals $ F,G $ on $ \cM(A) $, 
we say that 
\begin{align}
F(f) \asymp G(f) 
\quad \text{for $ f \in \cM(A) $} 
\label{}
\end{align}
if there exist positive constants $ c_1,c_2 $ such that 
\begin{align}
c_1 F(f) \le G(f) \le c_2 F(f) 
\quad \text{for $ f \in \cM(A) $} . 
\label{}
\end{align}

Now let us present several examples of norms which satisfy {\rm (A1)}-{\rm (A4)}. 

\begin{Ex}[Max norms]
Consider $ \| x \|^{(d)}_{\infty } = \max_i |x^i| $. 
It is obvious that the conditions {\rm (A1)}-{\rm (A3)} are satisfied. 
In fact, the partition of $ \partial B(0;1) $ in {\rm (A3)} 
can be obtained by separating $ \partial B(0;1) $ by hyperplanes 
$ \{ x \in \bR^d : x^i=j/k \} $ 
for $ i=1,\ldots,d $ and $ j=-k,\ldots,k $. 
Let us study $ N(k)= N^{(d)}_{\infty }(k) $ and its asymptotic behavior. 
For $ k \in \bZ_{\ge 1} $, we have 
\begin{align}
N^{(d)}_{\infty }(k) =& 
\sharp \{ x \in \bZ^d : \| x \| \le k \} 
- 
\sharp \{ x \in \bZ^d : \| x \| \le k-1 \} 
\label{} \\
=& (2k+1)^d - (2k-1)^d . 
\label{eq: max norm}
\end{align}
Now we obtain 
\begin{align}
N^{(d)}_{\infty }(k) \sim d 2^d k^{d-1} 
\quad \text{as $ k \to \infty $} . 
\label{}
\end{align}
\end{Ex}

\begin{Ex}[$ \ell^1 $-norms]
Consider 
\begin{align}
\| x \|^{(d)}_1 = \sum_{i=1}^d |x^i| 
, \quad x \in \bR^d . 
\label{}
\end{align}
It is obvious that the conditions {\rm (A1)}-{\rm (A3)} are satisfied. 
In this case, 
\begin{align}
N(k) = N^{(d)}_1(k) = \sharp \{ x \in \bZ^d : \| x \|^{(d)}_1 = k \} 
\label{}
\end{align}
satisfies the recursive relation 
\begin{align}
N^{(d)}_1(k) = \sum_{j=0}^k N^{(1)}_1(j) N^{(d-1)}_1(k-j) 
, \quad d \ge 2 , \ k \ge 0 
\label{eq: recursive}
\end{align}
with initial condition 
\begin{align}
N^{(1)}_1(k) = 
\begin{cases}
1 \quad & \text{if $ k=0 $}, \\
2 \quad & \text{if $ k \ge 1 $} . 
\end{cases}
\label{eq: initial}
\end{align}
Since the moment generating function may be computed as 
\begin{align}
\sum_{k=0}^{\infty } s^k N^{(d)}_1(k) = \rbra{ \frac{1+s}{1-s} }^d 
, \quad 0<s<1 , 
\label{}
\end{align}
we see, by Tauberian theorem (see, e.g., \cite[Thm.XIII.5.5]{Feller}), that 
\begin{align}
N^{(d)}_1(k) 
\sim \frac{2^d}{(d-1)!} k^{d-1} 
\quad \text{as $ k \to \infty $}. 
\label{}
\end{align}
\end{Ex}

\begin{Ex}[Weighted $ \ell^1 $-norms]
Consider 
\begin{align}
\| x \|^{(d)}_{\rm w1} = \sum_{i=1}^d i |x^i| 
, \quad x \in \bR^d . 
\label{}
\end{align}
The conditions {\rm (A1)}-{\rm (A3)} are obviously satisfied. 

Now let us discuss the asymptotic behavior of $ N^{(d)}_{\rm w1}(k) $. 
Note that 
\begin{align}
N(k) = N^{(d)}_{\rm w1}(k) = \sharp \{ x \in \bZ^d : \| x \|^{(d)}_{\rm w1} = k \} 
\label{}
\end{align}
satisfies the recursive relation 
\begin{align}
N^{(d)}_{\rm w1}(k) 
= \sum_{j \in \bZ_{\ge 0}: \, k-dj \ge 0} N^{(1)}_1(j) N^{(d-1)}_{\rm w1}(k-dj) 
, \quad d \ge 2 , \ k \ge 0 
\label{eq: recursive2}
\end{align}
with initial condition $ N^{(1)}_{\rm w1}(k) \equiv N^{(1)}_1(k) $. 
Then, by induction, we can easily see that 
\begin{align}
| N^{(d)}_{\rm w1}(k) - a^{(d)} k^{d-1} | \le b^{(d)} k^{d-2} 
, \quad k \in \bZ_{\ge 1} , \ d \ge 2 
\label{}
\end{align}
for some positive constants $ a^{(d)},b^{(d)} $ 
where $ a^{(d)} $ is defined recursively as 
\begin{align}
a^{(1)} = 2 
, \quad 
a^{(d)} = \frac{2}{d(d-1)} a^{(d-1)} \ (d \ge 2) . 
\label{}
\end{align}
In particular, we see that 
$ N^{(d)}_{\rm w1}(k) \sim a^{(d)} k^{d-1} $ as $ k \to \infty $. 
For instance, by easy computations, we obtain 
\begin{align}
N^{(2)}_{\rm w1}(k) = 
\begin{cases}
1 \quad & \text{if $ k=0 $}, \\
2k \quad & \text{if $ k \ge 1 $} 
\end{cases}
\label{}
\end{align}
and 
\begin{align}
N^{(3)}_{\rm w1}(k) = 
\begin{cases}
1 \quad & \text{if $ k=0 $}, \\
\frac{2}{3} k^2 + 2 \quad & \text{if $ k \equiv 0 $ modulo 3, $ k \neq 0 $} , \\
\frac{2}{3} k^2 + \frac{4}{3} \quad & \text{if $ k \equiv 1,2 $ modulo 3} . 
\end{cases}
\label{}
\end{align}
\end{Ex}

\begin{Ex}[Transformation by unimodular matrices]
Let $ A $ be a unimodular $ d \times d $ matrix, 
i.e., $ A $ is a $ d \times d $ matrix 
whose entries are integers and whose determinant is 1 or $ -1 $. 
Let $ \| \cdot \|_0 $ be a norm on $ \bR^d $ satisfying {\rm (A1)}-{\rm (A3)}. 
Then the norm $ \| x \| = \| Ax \|_0 $ also satisfies {\rm (A1)}-{\rm (A3)}. 
Note that 
\begin{align}
\sharp \{ x \in \bZ^d : \| x \| = k \} 
= \sharp \{ x \in \bZ^d : \| x \|_0 = k \} 
, \quad k \in \bZ_{\ge 0} . 
\label{}
\end{align}
For example, the norm on $ \bR^3 $ defined as 
\begin{align}
\| (x^1,x^2,x^3) \| = |x^1-x^2| + |x^2-x^3| + |x^1-x^2+x^3| 
\label{}
\end{align}
satisfies {\rm (A1)}-{\rm (A3)}. 
\end{Ex}

\begin{Rem}
Let us consider the norm $ 2 \| x \|^{(d)}_{\infty } $. 
Then the conditions {\rm (A1)}-{\rm (A2)} are satisfied, 
but neither of {\rm (A3)} nor {\rm (A4)} is; in fact, 
\begin{align}
N(k) = 
\begin{cases}
N^{(d)}_{\infty }(k/2) \quad & \text{if $ k $ is even}, \\
0 & \text{if $ k $ is odd}. 
\end{cases}
\label{}
\end{align}
Nevertheless, we see that 
the conditions {\rm (I)}-{\rm (IV)} are equivalent to each other 
and also to 
\begin{align}
\sum_{k=1}^{\infty } k f(2k) < \infty , 
\label{}
\end{align}
which is strictly weaker than {\rm (V)} 
because there is no restriction on the values of $ f(2k+1) $. 
\end{Rem}

\section{Equivalence between {\rm (III)} and {\rm (IV)}} \label{sec: proof1}
\setcounter{equation}{0}

Let us introduce the {\em random walk local times}: 
\begin{align}
L^S_n(x) =& \sharp \cbra{ m = 1,2,\ldots,n : S_m = x } 
, \quad x \in \bZ^d , 
\label{} \\
L^{\| S \|}_n(k) =& \sharp \cbra{ m = 1,2,\ldots,n : \| S_m \| = k } 
, \quad k \in \bZ_{\ge 0} . 
\label{}
\end{align}
Then, for any non-negative function $ g $ on $ \bZ^d $, we have 
\begin{align}
\sum_{n=1}^{\infty } g(S_n) = \sum_{x \in \bZ^d} g(x) L^S_{\infty }(x) . 
\label{}
\end{align}
Taking the expectations of both sides, we have 
\begin{align}
E \sbra{ \sum_{n=1}^{\infty } g(S_n) } 
= \sum_{x \in \bZ^d} g(x) E \sbra{ L^S_{\infty }(x) } . 
\label{eq: E sum gS_n}
\end{align}
It is obvious by definition that 
\begin{align}
E \sbra{ L^S_{\infty }(x) } = \sum_{n=1}^{\infty } P(S_n=x) = G(0,x) 
\label{}
\end{align}
where $ G(x,y) $ is the {\em Green function} given as 
\begin{align}
G(x,y) = \sum_{n=1}^{\infty } P_x(S_n=y) . 
\label{}
\end{align}
Let $ | \cdot | $ denote the Euclidean norm of $ \bR^d $, i.e., 
$ |x|^2 = \sum_{i=1}^d (x^i)^2 $. 
We recall the following asymptotic behavior of the Green function: 

\begin{Thm}[{\cite{Spitzer}}] \label{thm: Spitzer}
It holds that 
\begin{align}
G(0,x) 
\sim \frac{\Gamma(d/2-1)}{2 \pi^{d/2}} |\det Q|^{-1/2} (x,Q^{-1}x)^{1-d/2} 
\quad \text{as $ |x| \to \infty $} . 
\label{}
\end{align}
In particular, if $ Q = \sigma^2 I $, then 
\begin{align}
|x|^{d-2} G(0,x) 
\to \frac{\Gamma(d/2-1)}{2 \pi^{d/2}} \sigma^{-2} 
\quad \text{as $ |x| \to \infty $} . 
\label{}
\end{align}
\end{Thm}

We can prove Theorem \ref{thm: Spitzer} 
in the same way as in Spitzer \cite[P26.1]{Spitzer}, 
so we omit the proof. 

\begin{Prop} \label{thm: asymp for Zd}
It holds that 
\begin{align}
E \sbra{ \sum_{n=1}^{\infty } g(S_n) } 
\asymp g(0) + \sum_{x \in \bZ^d \setminus \{ 0 \}} g(x) \| x \|^{2-d} 
\quad \text{for $ g \in \cM(\bZ^d) $} . 
\label{}
\end{align}
\end{Prop}

\begin{proof}
Since $ \| x \| \asymp |x| $ for $ x \in \bZ^d $, 
it follows from Theorem \ref{thm: Spitzer} that 
\begin{align}
G(0,x) \asymp \| x \|^{2-d} 
\quad \text{for $ x \in \bZ^d \setminus \{ 0 \} $}. 
\label{}
\end{align}
Combining it with \eqref{eq: E sum gS_n}, 
we obtain the desired result. 
\end{proof}

\begin{Rem}
It is now obvious from Proposition \ref{thm: asymp for Zd} that 
\begin{align}
\sum_{x \in \bZ^d} g(x) \| x \|^{2-d} < \infty 
\quad \text{implies} \quad 
P \rbra{ \sum_{n=1}^{\infty } g(S_n) < \infty } = 1 . 
\label{}
\end{align}
But we do not know whether the converse is true or not. 
\end{Rem}

The following proposition proves part of Theorem \ref{thm: main2}. 

\begin{Prop} \label{thm: III and IV}
Suppose that the condition {\rm (A1)} is satisfied. 
Then it holds that 
\begin{align}
E \sbra{ \sum_{n=1}^{\infty } f(\| S_n \|) } 
\asymp f(0) + \sum_{k=1}^{\infty } f(k) k^{2-d} N(k) 
\quad \text{for $ f \in \cM(\bZ_{\ge 0}) $} 
\label{}
\end{align}
and, in particular, that {\rm (III)} and {\rm (IV)} are equivalent. 
If, moreover, the condition {\rm (A4)} is satisfied, then it holds that 
\begin{align}
E \sbra{ \sum_{n=1}^{\infty } f(\| S_n \|) } 
\asymp f(0) + \sum_{k=1}^{\infty } k f(k) 
\quad \text{for $ f \in \cM(\bZ_{\ge 0}) $} 
\label{}
\end{align}
and, in particular, that {\rm (IV)} and {\rm (V)} are equivalent. 
\end{Prop}

\begin{proof}
The former half of Proposition \ref{thm: III and IV} 
is immediate from Propositions \ref{thm: asymp for Zd} and \ref{thm: main3} 
for $ g(x) = f(\| x \|) $. 
The latter half is immediate 
from Proposition \ref{thm: main3} below. 
\end{proof}

\begin{Rem}
Let $ p(x) $ denote the probability that the process visits $ x $ at least once: 
\begin{align}
p(x) = P(L^S_{\infty }(x) \ge 1) = P(T_x < \infty ) 
, \quad x \in \bZ^d 
\label{}
\end{align}
where $ T_x = \inf \{ n \ge 1 : S_n = x \} $ is the first hitting time of $ x $. 
Since $ L^S_{\infty }(x) = L^S_{\infty }(x) \circ \theta_{T_x} + 1 $ 
and by translation invariance, we may compute the distribution of the total local time 
$ L^S_{\infty }(x) $ as 
\begin{align}
P(L^S_{\infty }(x) \ge n) = p(x) p(0)^{n-1} 
, \quad x \in \bZ^d , \ n=1,2,\ldots 
\label{}
\end{align}
See \cite{MR} for some general discussions for symmetric Markov processes. 
Note that the Green function $ G(0,x) $ may be expressed as 
\begin{align}
G(0,x) = E \sbra{ L^S_{\infty }(x) } 
= \sum_{n=1}^{\infty } P(L^S_{\infty }(x) \ge n) 
= \frac{p(x)}{1-p(0)} . 
\label{}
\end{align}
\end{Rem}

\begin{Rem}
We do not know any explicit result 
about the law of the total local time $ L^{\| S \|}_{\infty }(k) $ 
for the norm process $ \| S \| $. 
\end{Rem}

\begin{Prop} \label{thm: main3}
Let $ \| \cdot \| $ be a norm on $ \bR^d $. 
Suppose that the condition {\rm (A4)} is satisfied. 
Then there exists $ k_1 \in \bZ_{\ge 1} $ such that 
$ N(k) \asymp k^{d-1} $ for $ k \ge k_1 $. 
\end{Prop}

\begin{proof}
By \eqref{eq: max norm}, we have 
\begin{align}
\sharp \cbra{ x \in \bZ^d : \| x \|^{(d)}_{\infty } \le k } 
\asymp k^d 
\quad \text{for $ k \in \bZ_{\ge 1} $} . 
\label{}
\end{align}
Note that $ \| x \| \asymp \| x \|^{(d)}_{\infty } $ for $ x \in \bZ^d $; 
in fact, any two norms on $ \bR^d $ are mutually equivalent. 
This immediately implies that 
\begin{align}
\sum_{j=0}^k N(j) = \sharp \cbra{ x \in \bZ^d : \| x \| \le k } 
\asymp k^d 
\quad \text{for $ k \in \bZ_{\ge 1} $} . 
\label{}
\end{align}
Hence there exist constants $ c_1,c_2 $ such that 
\begin{align}
c_1 k^d \le \sum_{j=0}^k N(j) \le c_2 k^d 
\quad \text{for $ k \in \bZ_{\ge 1} $} . 
\label{}
\end{align}

By the condition {\rm (A4)}, we have 
\begin{align}
k N(k) = \sum_{j=k+1}^{2k} N(k) \le \sum_{j=k+1}^{2k} N(j) 
\le c_2 (2k)^d 
\quad \text{for $ k \in \bZ_{\ge 1} $} . 
\label{}
\end{align}
Now we obtain $ N(k) \le c_3 k^{d-1} $ with $ c_3 = c_2 2^d $. 
Again by the condition {\rm (A4)}, we have 
\begin{align}
k N(k) = \sum_{j=1}^{k} N(k) \ge \sum_{j=0}^{k} N(j) 
\ge c_1 k^d 
\quad \text{for $ k \in \bZ_{\ge 1} $} . 
\label{}
\end{align}
Now we obtain $ N(k) \ge c_1 k^{d-1} $. 
This completes the proof. 
\end{proof}

\section{Proving that {\rm (I)} implies {\rm (IV)}} \label{sec: proof2}
\setcounter{equation}{0}

By the assumption {\rm (A2)}, 
we may identify each cell of $ B(0;r) $ with a subset of $ \bR^{d-1} $. 
So we may introduce the area measure $ \lambda $ on $ \partial B(0;1) $. 
We define $ \mu(\cdot) = \lambda(\cdot)/\lambda(\partial B(0;1)) $ 
and call it the {\em uniform measure} on $ \partial B(0;1) $. 

For $ k \in \bZ_{\ge 1} $, 
we define a probability measure on $ \bR^d $ by 
\begin{align}
\mu_k(A) 
= \frac{1}{N(k)} \sharp \cbra{ x \in k^{-1} \bZ^d \cap A : \| x \| = 1 } 
, \quad A \in \cB(\bR^d) . 
\label{}
\end{align}

\begin{Prop} \label{thm: muk conv weakly}
Suppose that {\rm (A1)}-{\rm (A3)} are satified. 
Then, as $ k \to \infty $, 
the measure $ \mu_k $ converges weakly to $ \mu $. 
\end{Prop}

\begin{proof}
Let $ \{ U^{(k)}_j:j=1,\ldots,M(k) \} $ be such as in the assumption {\rm (A3)}. 
Then we see that $ \mu(U^{(k)}_j) = M(k)^{-1} $ for any $ j $ and any $ k $. 
For $ j=1,\ldots,M(k) $, choose $ x^{(k)}_j \in U^{(k)}_j \cap (k^{-1} \bZ^d) $. 
We may choose $ \{ x^{(k)}_j:j=M(k)+1,\ldots,N(k) \} $ so that 
$ \{ x^{(k)}_j:j=1,\ldots,N(k) \} $ is an enumeration of the points 
of $ \{ x \in k^{-1} \bZ^d : \| x \| = 1 \} $. 

Let $ f:\bR^d \to \bR $ be a continuous function with compact support. 
It suffices to prove that 
\begin{align}
\int_{\bR^d} f(x) \mu_k(\d x) 
\stackrel{k \to \infty }{\longrightarrow} 
\int_{\partial B(0;1)} f(x) \mu(\d x) . 
\label{}
\end{align}
Note that 
\begin{align}
\int_{\bR^d} f(x) \mu_k(\d x) 
= \frac{1}{N(k)} \sum_{j=1}^{N(k)} f(x^{(k)}_j) . 
\label{}
\end{align}
Since $ M(k)/N(k) \to 1 $ as $ k \to \infty $, it suffices to prove that 
\begin{align}
\frac{1}{M(k)} \sum_{j=1}^{M(k)} f(x^{(k)}_j) 
\stackrel{k \to \infty }{\longrightarrow} 
\int_{\partial B(0;1)} f(x) \mu(\d x) . 
\label{}
\end{align}
Since $ \partial B(0;1) = \cup_j U^{(k)}_j $ 
and $ \mu(U^{(k)}_j) = M(k)^{-1} $, we obtain 
\begin{align}
& \absol{ \frac{1}{M(k)} \sum_{j=1}^{M(k)} f(x^{(k)}_j) 
- \int_{\partial B(0;1)} f(x) \mu(\d x) } 
\label{} \\
\le& \frac{1}{M(k)} \sum_{j=1}^{M(k)} 
\int_{U^{(k)}_j} \absol{ f(x^{(k)}_j) - f(x) } \mu(\d x) 
\label{} \\
\le& \max_{1 \le j \le M(k)} \max_{x,y \in U^{(k)}_j} \absol{ f(y) - f(x) } . 
\label{eq: last line}
\end{align}
By uniform continuity of $ f $ and by the assumption {\rm (A3)}, 
the quantity \eqref{eq: last line} converges to 0 as $ k \to \infty $. 
Therefore the proof is complete. 
\end{proof}

Let $ (B_t) $ denote 
a standard Brownian motion of dimension $ d \ge 3 $ starting from the origin. 
Set 
\begin{align}
g(x,y) = \int_0^{\infty } \frac{\d s}{(2 \pi s)^{d/2}} 
\exp \rbra{ -\frac{|x-y|^2}{2s} } 
, \quad x,y \in \bR^d . 
\label{}
\end{align}
For the uniform measure $ \mu $ on $ \partial B(0;1) $, we define 
\begin{align}
g \mu(x) = \int_{\bR^d} g(x,y) \mu (\d y) 
, \quad x \in \bR^d . 
\label{}
\end{align}
Then it is well-known (see \cite{Mey}; see also \cite[Thm.5.2.5]{FOT}) that 
there exists a unique positive continuous additive functional $ (L^{\mu}_t) $ such that 
\begin{align}
g \mu(\sigma B_t) - g \mu(\sigma B_0) + L^{\mu}_t 
\label{}
\end{align}
is a martingale with zero mean. 
The process $ (L^{\mu}_t) $ is called 
the {\em local time process on the union of cells $ \partial B(0;1) $ for $ (\sigma B_t) $}. 
The relation between the measure $ \mu $ 
and the positive continuous additive functional $ (L^{\mu}_t) $ is called 
the {\em Revuz correspondence} (see \cite{Rev}). 

The following theorem is an invariance principle for the random walk local time 
of the norm process. 

\begin{Thm} \label{thm: inv princ}
Suppose that {\rm (A0)}-{\rm (A3)} are satisfied. 
Then it holds that 
\begin{align}
\frac{L^{\| S \|}_{\infty }(k)}{k^{2-d} N(k)} 
\claw L^{\mu}_{\infty } 
\quad \text{as $ k \to \infty $} . 
\label{}
\end{align}
\end{Thm}

\begin{proof}
Note that 
\begin{align}
\frac{L^{\| S \|}_{\infty }(k)}{k^{2-d} N(k)} 
= k^{d-2} \sum_{n=1}^{\infty } \mu_k \rbra{ \cbra{ \frac{S_n}{k} } } . 
\label{}
\end{align}
Hence we obtain the desired result 
as an immediate consequence of Proposition \ref{thm: muk conv weakly} 
and Bass--Khoshnevisan \cite[Prop.6.3]{BK}. 
\end{proof}

Now we are in a position to prove 
that {\rm (I)} implies {\rm (IV)} in Theorem \ref{thm: main2}. 

\noindent
{\it Proof that {\rm (I)} implies {\rm (IV)} in Theorem \ref{thm: main2}.} 
Let us check that the assumptions of Proposition \ref{thm: limit Jeulin1} 
are satisfied for $ V(k) = L^{\| S \|}_{\infty }(k) $, 
$ \Phi(k) = k^{2-d} N(k) $ 
and $ X = L^{\mu}_{\infty } $. 

By Theorem \ref{thm: inv princ}, 
assumption \eqref{eq: conv in law} is satisfied. 

Let us show that $ P(L^{\mu}_{\infty } \le 0) = 0 $. 
The first hitting place of the union of cells $ \partial B(0;1) $ for the Brownian motion 
is almost surely an interior point of some cell of the $ d $-polytope $ B(0;1) $ 
by assumption {\rm (A2)}. 
Hence it holds that, starting afresh at the first hitting time, 
the local time on the union of cells $ \partial B(0;1) $ 
is locally equal to the local time on the hyperplane which contains the cell. 
Since the local time at the origin for one-dimensional Brownian motion 
is positive almost surely at any positive time, 
we see that $ L^{\mu}_{\infty } $ is positive almost surely. 

Thus we may apply Proposition \ref{thm: limit Jeulin1} (or Corollary \ref{thm: Shiga2}) 
and we see that {\rm (I)} implies {\rm (IV)}. 
The proof is now complete. 
\qed

\section{A remark on Jeulin's lemma} \label{sec: Shiga}
\setcounter{equation}{0}

The results of this section are mainly due to Tokuzo Shiga \cite{Shiga}. 

\subsection{Counterexample to Jeulin's lemma without $ E[X]<\infty $}

The following proposition gives a counterexample 
to Jeulin's lemma (Theorem \ref{thm: Jeulin}) without $ E[X]<\infty $. 

\begin{Prop}[\cite{Shiga}] \label{thm: Shiga5}
There exist a non-negative measurable process $ (X(t):0 < t \le 1) $, 
a positive function $ \varphi $ on $ (0,1] $, 
a random variable $ X $, 
and a non-negative Borel measure $ \mu $ on $ (0,1] $ 
such that 
\begin{align}
E[X]=\infty \quad \text{and} \quad P(X>0) = 1 , 
\label{eq: Shiga5-1}
\end{align}
\begin{align}
\frac{X(t)}{\varphi(t)} \law X 
\quad \text{holds for each fixed $ 0 < t \le 1 $}, 
\label{eq: Shiga5-2}
\end{align}
\begin{align}
\int_0^1 \varphi(t) \mu(\d t) < \infty 
\label{eq: Shiga5-3}
\end{align}
but 
\begin{align}
P \rbra{ \int_{\eps}^1 X(t) \mu(\d t) < \infty \ (\forall \eps>0) 
, \quad 
\int_0^1 X(t) \mu(\d t) = \infty } = 1 . 
\label{eq: Shiga5-4}
\end{align}
\end{Prop}

\begin{proof}
Let $ (X(t)) $ be an $ \alpha $-stable subordinator with $ 0<\alpha \le 1/2 $. 
Then we have \eqref{eq: Shiga5-1} and \eqref{eq: Shiga5-2} 
for $ \varphi(t) \equiv t^{1/\alpha } $. 
Set 
\begin{align}
\mu(\d t) = t^{-1-1/\alpha } (\log 1/t)^{-1/\alpha } \d t 
\label{}
\end{align}
so that $ \mu((t,1])^{\alpha } \sim C t^{-1} (\log 1/t)^{-1} $ 
as $ t \to 0+ $ for some positive constant $ C $. 
Thus we obtain \eqref{eq: Shiga5-3}. 
Since we have 
\begin{align}
E \sbra{ \exp - \int_0^1 X(t) \mu(\d t) } = \exp - \int_0^1 \mu((t,1])^{\alpha } \d t = 0 , 
\label{}
\end{align}
we obtain \eqref{eq: Shiga5-4}. 
\end{proof}

\subsection{A limit version of Jeulin's lemma}

\begin{Thm}[\cite{Shiga}] \label{thm: Shiga}
Let $ (X(t):0 < t \le 1) $ be a non-negative measurable process, 
$ \varphi $ a positive function defined on $ (0,1] $, 
and $ \mu $ a non-negative Borel measure on $ (0,1] $. 
Suppose that there exists a random variable $ X $ 
with $ P(X>0)>0 $ such that 
\begin{align}
\frac{X(t)}{\varphi(t)} \claw X 
\quad \text{as} \ t \to 0+ . 
\label{eq: Shiga3}
\end{align}
Suppose, moreover, that 
\begin{align}
\int_{\eps}^1 \varphi(t) \mu(\d t) < \infty 
\quad \text{for every $ 0<\eps<1 $}. 
\label{eq: Shiga2}
\end{align}
Then it holds that 
\begin{align}
P \rbra{ \int_0^1 X(t) \mu(\d t) < \infty } = 1 
\quad \text{implies} \quad 
\int_0^1 \varphi(t) \mu(\d t) < \infty . 
\label{eq: Shiga1}
\end{align}
\end{Thm}

\begin{proof}
Suppose that 
\begin{align}
P \rbra{ \int_0^1 X(t) \mu(\d t) < \infty } = 1 
\label{eq: Shiga4}
\end{align}
but that $ \int_0^1 \varphi(s) \mu(\d s) = \infty $. 
For each $ \eps>0 $, we define a probability measure $ \mu_{\eps} $ by 
\begin{align}
\mu_{\eps}(\d t) = C_{\eps}^{-1} 1_{(\eps,1]}(t) \varphi(t) \mu(\d t) 
\quad \text{with} \quad 
C_{\eps} = \int_{\eps}^1 \varphi(t) \mu(\d t) 
\label{}
\end{align}
where $ C_{\eps} $ is finite by the assumption \eqref{eq: Shiga2}. 
Then $ C_{\eps} \to \infty $ and $ \mu_{\eps} \cdist \delta_0 $ as $ \eps \to 0+ $, 
where $ \delta_0 $ stands for the unit point mass at $ 0 $. 
Using Jensen's inequality and changing the order of integration, we have 
\begin{align}
E \sbra{ \exp - C_{\eps}^{-1} \int_{\eps}^1 X(t) \mu(\d t) } 
=& 
E \sbra{ \exp - \int_{\eps}^1 \frac{X(t)}{\varphi(t)} \mu_{\eps}(\d t) } 
\label{} \\
\le& 
\int_{\eps}^1 E \sbra{ \exp - \frac{X(t)}{\varphi(t)} } \mu_{\eps}(\d t) . 
\label{}
\end{align}
Hence it follows from \eqref{eq: Shiga4} and \eqref{eq: Shiga3} that 
\begin{align}
1 \le \lim_{t \to 0+} E \sbra{ \exp - \frac{X(t)}{\varphi(t)} } 
= E \sbra{ \e^{-X} } , 
\label{}
\end{align}
which implies $ P(X=0)=1 $. 
This is a contradiction to the assumption that $ P(X>0)>0 $. 
\end{proof}

From Theorem \ref{thm: Shiga}, 
we obtain another version of Jeulin's lemma in discrete time. 

\begin{Cor} \label{thm: Shiga2}
Let $ (V(k):k \in \bZ_{\ge 1}) $ be a non-negative measurable process 
and $ \Phi $ a positive function on $ \bZ_{\ge 1} $. 
Suppose that 
there exists a random variable $ X $ with 
\begin{align}
P(X>0)>0 
\label{}
\end{align}
such that 
\begin{align}
\frac{V(k)}{\Phi(k)} \claw X 
\quad \text{as $ k \to \infty $}. 
\label{}
\end{align}
Then, for any non-negative function $ f $ on $ \bZ_{\ge 1} $, 
it holds that 
\begin{align}
P \rbra{ \sum_{k=1}^{\infty } f(k) V(k) < \infty } = 1 
\quad \text{implies} \quad 
\sum_{k=1}^{\infty } f(k) \Phi(k) < \infty. 
\label{eq: Shiga Jeulin}
\end{align}
\end{Cor}

\begin{proof}
Take 
\begin{align}
X(t) = V([1/t]) 
, \quad 
\varphi(t) = \Phi([1/t]) 
\label{}
\end{align}
where $ [x] $ stands for the smallest integer which does not exceed $ x $ 
and 
\begin{align}
\mu = \sum_{k=1}^{\infty } f(k) \delta_{1/k} . 
\label{}
\end{align}
Then the desired result is immediate from Theorem \ref{thm: Shiga}. 
\end{proof}

Proposition \ref{thm: limit Jeulin1} and Corollary \ref{thm: Shiga2} 
cannot be unified in the following sense: 

\begin{Prop}
There exist a non-negative measurable process $ (V(k):k \in \bZ_{\ge 1}) $, 
a positive function $ \Phi $ on $ \bZ_{\ge 1} $, 
a random variable $ X $, 
and a non-negative function $ f $ on $ \bZ_{\ge 1} $ 
such that 
\begin{align}
P(X>0)>0 , \quad 
\frac{V(k)}{\Phi(k)} \claw X 
\ \text{as $ k \to \infty $} , 
\label{eq: counterex1}
\end{align}
and 
\begin{align}
P \rbra{ \sum_{k=1}^{\infty } f(k) V(k) < \infty } > 0 
\label{}
\end{align}
but 
\begin{align}
\sum_{k=1}^{\infty } f(k) \Phi(k) = \infty. 
\label{eq: counterex3}
\end{align}
\end{Prop}

\begin{proof}
Let $ X $ be such that 
\begin{align}
P(X=0) = P(X=1) = \frac{1}{2} 
\label{}
\end{align}
and set $ V(k) = X $ for $ k \in \bZ_{\ge 1} $. 
Then we have \eqref{eq: counterex1}-\eqref{eq: counterex3} 
for $ \Phi(k) \equiv 1 $ and $ f(k) \equiv 1 $. 
\end{proof}

\subsection{A counterexample}

We give a counterexample 
to the converse of \eqref{eq: Shiga Jeulin} 
where the assumptions of Corollary \ref{thm: Shiga2} are satisfied. 

\begin{Prop}[\cite{Shiga}] \label{thm: Shiga3}
There exist a non-negative measurable process $ (V(k):k \in \bZ_{\ge 1}) $, 
a positive function $ \Phi $ on $ \bZ_{\ge 1} $, 
and a non-negative function $ f $ on $ \bZ_{\ge 1} $ 
such that 
\begin{align}
\frac{V(k)}{\Phi(k)} \claw 1 
\quad \text{as $ k \to \infty $} 
\label{eq: counterex4}
\end{align}
and 
\begin{align}
\sum_{k=1}^{\infty } f(k) \Phi(k) < \infty. 
\label{eq: counterex5}
\end{align}
but 
\begin{align}
P \rbra{ \sum_{k=1}^{\infty } f(k) V(k) = \infty } = 1 . 
\label{eq: counterex6}
\end{align}
\end{Prop}

\begin{proof}
Let $ 0<\alpha <1/2 $. 
Let $ (V_0(k)) $ be a sequence of i.i.d.~random variables such that 
\begin{align}
E[\e^{- \lambda V_0(k)}] = \e^{- \lambda^{\alpha }} 
, \quad \lambda>0 , \ k \in \bZ_{\ge 1} . 
\label{}
\end{align}
Set $ \Phi(k) \equiv k $ and $ f(k) \equiv k^{-1/\alpha } $. 
Then \eqref{eq: counterex5} holds and we have 
\begin{align}
\frac{V_0(k)}{\Phi(k)} \claw 0 
\quad \text{as $ k \to \infty $} . 
\label{}
\end{align}
Since we have 
\begin{align}
E \sbra{ \exp - \sum_{k=1}^{\infty } f(k) V_0(k) } 
= \prod_{k=1}^{\infty } E \sbra{ \e^{- f(k) V_0(k)} } 
= \exp - \sum_{k=1}^{\infty } k^{-1} = 0 , 
\label{}
\end{align}
we obtain 
\begin{align}
P \rbra{ \sum_{k=1}^{\infty } f(k) V_0(k) = \infty } = 1 . 
\label{}
\end{align}
Since we may take $ V(k) = k + V_0(k) $, 
we also obtain \eqref{eq: counterex4} and \eqref{eq: counterex6}. 
The proof is now complete. 
\end{proof}

{\bf Acknowledgements.} 
The authors would like to thank Professor Tokuzo Shiga 
who kindly allowed them to append to this paper 
his detailed study \cite{Shiga} about Jeulin's lemma. 
They also thank Professors Marc Yor, Katsushi Fukuyama and Patrick J. Fitzsimmons 
for valuable comments. 
They are thankful to the referee for pointing out several errors in the earlier version. 
The first author, Ayako Matsumoto, expresses her sincerest gratitudes 
to Professors Yasunari Higuchi and Taizo Chiyonobu 
for their encouraging guidance in her study of mathematics. 
The research of the second author, Kouji Yano, 
was supported by KAKENHI (20740060).

\end{document}